\documentclass{amsart}

\usepackage{amsthm,amsmath,amssymb}
\usepackage{mathabx}
\usepackage[margin=1in]{geometry}
\usepackage{datetime}
\usepackage{graphicx}
\usepackage[numbers]{natbib}
\usepackage{booktabs}
\usepackage{longtable}
\usepackage{float}

\newcommand\Tstrut{\makebox[0pt][c]{\rule{0pt}{2.6ex}}}		
\newcommand\Bstrut{\makebox[0pt][c]{\rule[-1.2ex]{0pt}{0pt}}}	

\makeatletter
\patchcmd{\@settitle}{\uppercasenonmath\@title}{}{}{}
\patchcmd{\@setauthors}{\MakeUppercase}{}{}{}
\patchcmd{\section}{\scshape}{}{}{}
\makeatother

\begin{document}

\title{Some new parameterizations for the Diophantine bi-orthogonal monoclinic piped}
\author{Randall L. Rathbun}
\email{randallrathbun@gmail.com}

\begin{abstract}
The bi-orthogonal monoclinic Diophantine parallelepiped is introduced, then the $s$-parameters and their
governing equation for the bi-orthogonal monoclinic Diophantine parallelepiped are discussed.
Previous discoveries and parameterizations of the monoclinic piped are noted.
Then two parameterizations $P\left[\frac{1}{2},s_2,s_3,s_4\right], s_i \in \mathbb{Q},\mathbb{Z}$ are given
for a specific type of Diophantine bi-orthogonal monoclinic parallelepiped.
Next, a parameterization $P\left[s_1,s_2,s_3,s_4\right], s_i \in \mathbb{Q},\mathbb{Z}$ is presented
which covers 99.5\% of solutions found by raw computer searches. Several asymptotic sequences approaching
the {\it perfect cuboid} are listed, and some final comments made.
\end{abstract}

\maketitle

\setcounter{section}{0}
\section*{{\bf The bi-orthogonal monoclinic Diophantine parallelepiped}}

\begin{figure}[!ht]
\centering
\includegraphics[scale=1.0]{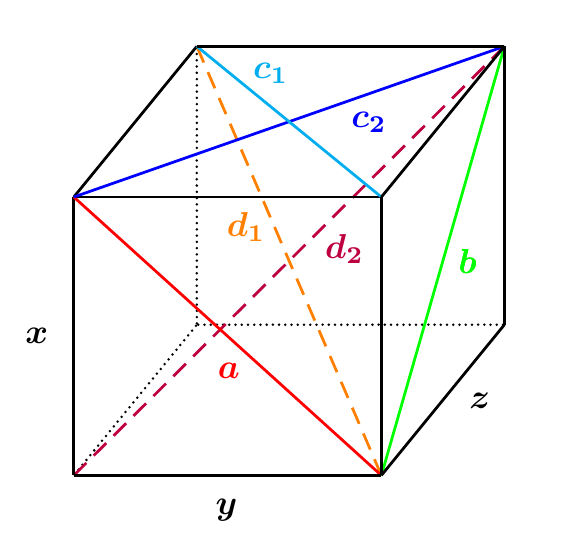}
\vspace{-6pt}
\caption{The Bi-orthogonal Monoclinic Cuboid.}
\label{fig:mono}
\end{figure}

Let Figure \ref{fig:mono} be a general bi-orthogonal monoclinic Diophantine parallelepiped, which is a
cuboid with two right angles at any vertice(bi-orthogonal), but the third is not a right angle,
but either acute or obtuse(monoclinic). This figure is composed of two congruent parallelograms (not
rectangles) joined by 4 orthogonal rectangles, all of whose edges are rational or integer.
Let $x$, $y$, $z$, denote the three different edges, with $x,y,z \in Z$.
The face rectangle $(x,y)$ has diagonal $a$, the face rectangle $(x,z)$ has diagonal $b$
and the face parallelogram $(y,z)$ has diagonals $c_1, c_2$.
The two different body diagonals are denoted by $d_1, d_2$.

These lengths satisfy the equations
\begin{align}
x^2 + y^2 & = a^2 \\
x^2 + z^2 & = b^2 \\
x^2 + c_1^2 & = d_1^2 \\
x^2 + c_2^2 & = d_2^2 \\
2y^2 + 2z^2 & = c_1^2 + c_2^2 \\
2y^2 + 2b^2 & = d_1^2 + d_2^2 \\
2a^2 + 2z^2 & = d_1^2 + d_2^2
\end{align}

These equations are from Wyss \cite{wyss3}.

The bi-orthogonal monoclinic Diophantine parallelepiped belongs to the family of {\it perfect parallelepipeds}. These
perfect parallelepipeds have all twelve edges rational, (there are only 3 distinct values), all six face diagonals rational,
and all four body diagonals rational. The family contains triclinic, biclinic, and monoclinic pipeds. The {\it perfect cuboid},
if it exists, would belong to this family. It is still unknown if this cuboid exists or not, see Guy \cite{guy}.

\section*{{\bf The $s$-parameters and the governing equation}}

We divide the integer or rational lengths of the piped by $x$ to obtain:
\begin{align}
1 + \left(\frac{y}{x}\right)^2 & = \left(\frac{a}{x}\right)^2
& \text{ let } \frac{y}{x} = u_1 \text{ \& } \frac{a}{x} = v_1 \text{ then } 1 + u_1^2 = v_1^2 \\
1 + \left(\frac{z}{x}\right)^2 & = \left(\frac{b}{x}\right)^2
& \text{ let } \frac{z}{x} = u_2 \text{ \& } \frac{b}{x} = v_2 \text{ then } 1 + u_2^2 = v_2^2 \\
1 + \left(\frac{c_1}{x}\right)^2 & = \left(\frac{d_1}{x}\right)^2
& \text{ let } \frac{c_1}{x} = u_3 \text{ \& } \frac{d_1}{x} = v_3 \text{ then } 1 + u_3^2 = v_3^2 \\
1 + \left(\frac{c_2}{x}\right)^2 & = \left(\frac{d_2}{x}\right)^2
& \text{ let } \frac{c_2}{x} = u_4 \text{ \& } \frac{d_2}{x} = v_4 \text{ then } 1 + u_4^2 = v_4^2
\end{align}
We have four rational equations for Pythagorean triangles.
\begin{equation}
1 + u_i^2 = v_i^2 , \quad i=1\dots 4
\end{equation}
which can be solved in parametric form:
\begin{equation}
u_k = \frac{1-s_k^2}{2s_k}, \quad v_k = \frac{1+s_k^2}{2s_k} \quad \text{ for } k=1 \dots 4
\end{equation}
for $s\in \mathbb{Q}$. We call the quadruple $s_{i=1\dots 4}$, the $s$-parameters $ = \left[s_1,s_2,s_3,s_4\right]$.
Wyss provides the derivation for equations (8-13) in his paper \cite{wyss3}.

If we substitute the $s$-parameters into equations (5-7), we obtain:
\begin{align}
2u_1^2 + 2u_2^2 & = u_3^2 + u_4^2 \\
2u_1^2 + 2v_2^2 & = v_3^2 + v_4^2 \\
2v_1^2 + 2u_2^2 & = v_3^2 + v_4^2
\end{align}
We next substitute $s_{i=1\dots 4}$ for the $u_{i=1\dots 4}$ parameters into equation (14) to derive the governing equation:
\begin{equation}
2 \left(\frac{1-s_1^2}{2s_1}\right)^2 + 2\left(\frac{1-s_2^2}{2s_2}\right)^2 = \left(\frac{1-s_3^2}{2s_3}\right)^2 + \left(\frac{1-s_4^2}{2s^4}\right)^2
\end{equation}
and simplify equation (17) to obtain:
\begin{align}
\label{eq:govern}
\begin{split}
& s_1^2s_2^2s_3^4s_4^2 + s_1^2s_2^2s_3^2s_4^4 - 2s_1^4s_2^2s_3^2s_4^2 - 2s_1^2s_2^4s_3^2s_4^2 + \\
& + 4s_1^2s_2^2s_3^2s_4^2 - 2s_1^2s_3^2s_4^2 - 2s_2^2s_3^2s_4^2 + s_1^2s_2^2s_3^2 + s_1^2s_2^2s_4^2 = 0
\end{split}
\end{align}
which is the governing equation for the $s$-parameters. This equation (\ref{eq:govern}) is given by Sharipov in \cite{shar}.

Ruslan Sharipov \cite{shar} points out that these parameters have to satisfy the following ten polynomial inequalities:
\begin{align}
\label{eq:satisfy}
\begin{split}
& s_1s_2^2s^3 + s_1^2s_2s_3 - s_1s_2s_3^2 + s_1s_2 - s_2s_3 - s_1s_3 < 0, \\
& s_1s_2^2s^4 + s_1^2s_2s_4 - s_1s_2s_4^2 + s_1s_2 - s_2s_4 - s_1s_4 < 0. \\
& 0 < s_i < 1 \; \text{ for } i=1,2,3,4
\end{split}
\end{align}
for rational Diophantine monoclinic pipeds to exist.

Thus we conclude: {\it Every rational Diophantine monoclinic piped corresponds to some
 $\left[s_1,s_2,s_3,s_4\right]$ s-parameter satisfying equations (\ref{eq:govern},\ref{eq:satisfy})}.

\section*{{\bf Previous parameterizations of the perfect parallelepiped}}

About 2009, Jorge Sawyer and Clifford Reiter discovered a perfect parallelepiped \cite{saw}. Clifford Reiter and Jordan Tirrell
then released an additional series of papers on parameterizing them \cite{reit1,reit2}. A more recent paper by Sokolowshy,
VanHooft, Volkert, and Reiter has also appeared, which provides an infinite family of monoclinic Diophantine
pipeds \cite{soko}.

Walter Wyss was able to successfully parameterize the perfect parallelogram which has rational sides and two
rational diagonals \cite{wyss1,wyss2}, and he extended this to parameterize a family of bi-orthogonal monoclinic
pipeds \cite{wyss3}.

Ruslan Sharipov recast Wyss' parametrizations, as equations (6.24-6.27), and extended them to three additional mappings for
the bi-orthogonal monoclinic piped, in his equation (6.31), see Sharipov \cite{shar}.

Based upon Wyss's and Sharipov's works, the author recently discovered four parametrizations of the monoclinic piped \cite{rath}.

\section*{{\bf Two parameterizations of the specific $s$-parameter $\left[\frac{1}{2},s_2,s_3,s_4\right]$ }}

This author began a computer search for $s$-parameter solutions to equation (\ref{eq:govern}), which originally was
an $O^8$ type of search, very slow and inefficient. By studying the equation, some ways were found to drop the search space
to an $0^6$ or even an $O^4$ type of search which was much more efficient in discovering solutions.

In particular, it was decided to set $s_1=\frac{1}{2}$ and look for solutions.

After finding about a hundred or so solutions to $s=\left[\frac{1}{2},s_2,s_3,s_4\right]$, it was noticed that the solutions
seemed to occur in two patterns.

{\bf First pattern type $s=\left[\frac{1}{2},\frac{c}{a},\frac{c}{b},\frac{a}{b}\right]$ }

Substituting these $s$-parameters into equation (\ref{eq:govern}), we obtain:
\begin{align}
\frac{a^2 - 2b^2}{4b^6a^2} c^6 + \frac{2a^4 - 9b^2a^2 + 2b^4}{8b^6a^2}c^4 + \left(\frac{1}{4b^2} - \frac{1}{2b^4} a^2 \right)c^2 & = 0 \notag \\
\text{or} \qquad \qquad \qquad \qquad \qquad \notag \\
\label{eq:pat1}
\frac{c^2\left(2c^2-b^2+2a^2\right)\left(2b^2c^2-a^2c^2+2a^2b^2\right)}{8a^2b^6} & = 0
\end{align}
which has 3 solutions for $c$, found by setting each of the 3 factors of the numerator of equation(\ref{eq:pat1}) $=0$.

First solution for the first factor, $c^2 = 0 \;$:
\begin{equation}
	c = 0
\end{equation}
This solution is trivial, and thus ignored.

Second solution for the second factor, $2c^2-b^2 + 2a^2 = 0 \;$:
\begin{align}
2c^2 - b^2 + 2a^2 & = 0 \notag \\
2c^2 & = b^2 - 2a^2 \notag \\
4c^2 & = 2b^2 - 4a^2 \notag \\
c & = \pm \frac{1}{2}\sqrt{2b^2-4a^2}
\end{align}
We want $2b^2 - 4a^2 = \square$ which has a general parametric solution in $m,n$.
\begin{align*}
\alpha^2 + \beta^2 & = 2\gamma^2 \\
\alpha & = 2mn+m^2-n^2 \\
\beta & = 2mn+n^2-m^2 \\
\gamma & = m^2 + n^2
\end{align*}
Adapting that parametrization, we obtain the following expression for the second factor
$\left(2c^2-b^2+2a^2\right) = 0$ as:
\begin{align}
& a = mn+\frac{1}{2}\left(m^2-n^2\right) \qquad
b = m^2 + n^2 \qquad
c = mn+\frac{1}{2}\left(n^2-m^2\right) \notag \\
& \text{\hspace{2.0in} or} \label{eq:second} \\
& a = mn+\frac{1}{2}\left(n^2-m^2\right) \qquad
b = m^2 + n^2 \qquad
c = mn+\frac{1}{2}\left(m^2-n^2\right) \notag
\end{align}
Notice the swap between $m,n$ for $a,c$.

Third solution for the third factor, $2b^2c^2 - a^2c^2 + 2a^2b^2 = 0 \;$:
\begin{align}
2b^2c^2 - a^2c^2 + 2a^2b^2 & = 0 \notag \\
c^2(2b^2-a^2) & = -2a^2b^2 \notag \\
c^2(a^2-2b^2) & = 2a^2b^2 \notag \\
c^2 & = \frac{2a^2+b^2}{a^2-2b^2} = \frac{4a^2b^2}{2a^2-4b^2} \notag \\
c & = \pm \frac{2ab}{\sqrt{2a^2-4b^2}}
\end{align}
Notice the similarity in the denominator to the second solution for $c$, except $a,b$ are swapped.

Using the same parametric solution for $2a^2-4b^2 = \square$, we find for the third factor, the
expression:
\begin{align}
& a = m^2 + n^2 \qquad
b = mn+\frac{1}{2}\left(m^2-n^2\right) \qquad
c = \frac{m^4+2nm^3+2n^3m-n^4}{-m^2+2nm+n^2} \notag \\
& \text{\hspace{2.0in} or} \label{eq:third} \\
& a = m^2 + n^2 \qquad
b = mn+\frac{1}{2}\left(n^2-m^2\right) \qquad
c = \frac{-m^4+2nm^3+2n^3m+n^4}{m^2+2nm-n^2} \notag
\end{align}
Collecting the four parametrized solutions, we obtain the $s$-parameter sets as:
\begin{align}
		a_1 & = 2mn+m^2-n^2 \notag \\
		b_1 & = 2(m^2+n^2) \notag \\
		c_1 & = 2mn+n^2-m^2 \notag \\
		s_1(m,n) & = \left[\frac{1}{2},\frac{2mn+n^2-m^2}{2mn+m^2-n^2},\frac{2mn+n^2-m^2}{2(m^2+n^2)},\frac{2mn+m^2-n^2}{2(m^2+n^2)}\right]
\end{align}
\begin{align}
		a_2 & = 2mn+n^2-m^2 \notag \\
		b_2 & = 2(m^2+n^2) \notag \\
		c_2 & = 2mn+m^2-n^2 \notag \\
		s_2(m,n) & = \left[\frac{1}{2},\frac{2mn+m^2-n^2}{2mn+n^2-m^2},\frac{2mn+m^2-n^2}{2(m^2+n^2)},\frac{2mn+n^2-m^2}{2(m^2+n^2)}\right]
\end{align}
\begin{align}
		a_3 & = 2(n^2+2nm-m^2)(m^2+n^2) \notag \\
		b_3 & = (n^2+2nm-m^2)(2mn+m^2-n^2) \notag \\
		c_3 & = 2(m^4+2nm^3+2n^3m-n^4) \notag \\
		s_3(m,n) & = \left[\frac{1}{2},\frac{m^2+2mn-n^2}{n^2+2nm-m^2},\frac{2(m^2+n^2)}{n^2+2nm-m^2},\frac{2(m^2+n^2)}{m^2+2mn-n^2}\right]
\end{align}
\begin{align}
		a_4 & = 2(m^2+2nm-n^2)(m^2+n^2) \notag \\
		b_4 & = (m^2+2nm-n^2)(2mn+n^2-m^2) \notag \\
		c_4 & = 2(-m^4+2nm^3+2n^3m+n^4) \notag \\
		s_4(m,n) & = \left[\frac{1}{2},\frac{m^2-2mn-n^2}{n^2-2nm-m^2},\frac{2(m^2+n^2)}{m^2+2mn-n^2},\frac{2(m^2+n^2)}{n^2+2nm-m^2}\right]
\end{align}
The parametric solutions are permutations of $s$-parameter sets, so only 1 is needed.
For example, using $(m,n)=(2,1)$:
\begin{align*}
s_1(2,1) & = [1/2,1/7,1/10,7/10] \\
s_2(2,1) & = [1/2,7,7/10,1/10] \\
s_3(2,1) & = [1/2,7,10,10/7] \\
s_4(2,1) & = [1/2,1/7,10/7,10]
\end{align*}
We choose $s_1(m,n)$ as the best parametrization for our purposes.

In summary, for the first pattern type, $s=\left[\frac{1}{2},\frac{c}{a},\frac{c}{b},\frac{a}{b}\right]$,
we can choose for integers $m,n \in \mathbb{Z}$, the following rational solution for the $s$-parameter:
\begin{equation}
	s(m,n) = \left[\frac{1}{2},\frac{2mn+n^2-m^2}{2mn+m^2-n^2},\frac{2mn+n^2-m^2}{2(m^2+n^2)},\frac{2mn+m^2-n^2}{2(m^2+n^2)}\right]
\label{eq:sintg}
\end{equation}
For a rational parameter $q \in \mathbb{Q}$, we substitute into equation (\ref{eq:sintg}), $m=q$, $n=1$ and obtain:
\begin{equation}
	s(q) = \left[\frac{1}{2}, \frac{-q^2 + 2q + 1}{q^2 + 2q - 1}, \frac{-q^2 + 2q + 1}{2q^2 + 2}, \frac{q^2 + 2q - 1}{2q^2 + 2}\right]
\label{eq:srat}
\end{equation}
The results from $s(q)$ match the $s$-parameters recoverable from the parameterization of equations (8-11) in this author's paper \cite{rath}, after using the $x$, $y$, $z$, $c_1$, and $c_2$ values to recover the corresponding $s_i$ term.

{\bf Second pattern type $s=\left[\frac{1}{2},\frac{d}{b},\frac{d}{a},\frac{d}{c}\right]$ }

A laborious computer search found the following solutions for $s_1=\frac{1}{2}$, in table \ref{table:search}. which match the second pattern.
\begin{longtable}{lllc}
\hspace{0.5cm} s-series & $b,a,c$ & $d$ & $q$ \\
\toprule
$\left[\frac{1}{2},\frac{16}{7},\frac{16}{5},\frac{16}{35}\right]$ & 7,5,35 & 16 & $\frac{1}{3}$ \Tstrut\Bstrut \\
$\left[\frac{1}{2},\frac{80}{119},\frac{80}{91},\frac{80}{221}\right]$ & 119,91,221 & 80 & $\frac{1}{5}$ \Tstrut\Bstrut \\
$\left[\frac{1}{2},\frac{160}{161},\frac{160}{119},\frac{160}{391}\right]$ & 161,119,391 & 160 & $\frac{1}{4}$ \Tstrut\Bstrut \\
$\left[\frac{1}{2},\frac{224}{527},\frac{224}{425},\frac{224}{775}\right]$ & 527,425,775 & 224 & $\frac{1}{7}$ \Tstrut\Bstrut \\
$\left[\frac{1}{2},\frac{480}{1519},\frac{480}{1271},\frac{480}{2009}\right]$ & 1519,1271,2009 & 480 & $\frac{1}{9}$ \Tstrut\Bstrut \\
$\left[\frac{1}{2},\frac{560}{41},\frac{560}{29},\frac{560}{1189}\right]$ & 41,29,1189 & 560 & $\frac{2}{5}$ \Tstrut\Bstrut \\
$\left[\frac{1}{2},\frac{560}{1081},\frac{560}{851},\frac{560}{1739}\right]$ & 1081,851,1739 & 560 & $\frac{1}{6}$ \Tstrut\Bstrut \\
$\left[\frac{1}{2},\frac{880}{3479},\frac{880}{2989},\frac{880}{4331}\right]$ & 3479,2989,4331 & 880 & $\frac{1}{11}$ \Tstrut\Bstrut \\
$\left[\frac{1}{2},\frac{1344}{3713},\frac{1344}{3055},\frac{1344}{5135}\right]$ & 3713,3055,5135 & 1344 & $\frac{1}{8}$ \Tstrut\Bstrut \\
$\left[\frac{1}{2},\frac{1456}{6887},\frac{1456}{6035},\frac{1456}{8245}\right]$ & 6887,6035,8245 & 1456 & $\frac{1}{13}$ \Tstrut\Bstrut \\
$\left[\frac{1}{2},\frac{1680}{1241},\frac{1680}{901},\frac{1680}{3869}\right]$ & 1241,901,3869 & 1680 & $\frac{2}{7}$ \Tstrut\Bstrut \\
$\left[\frac{1}{2},\frac{2240}{12319},\frac{2240}{10961},\frac{2240}{14351}\right]$ & 12319,10961,14351 & 2240 & $\frac{1}{15}$ \Tstrut\Bstrut \\
$\left[\frac{1}{2},\frac{2464}{2047},\frac{2464}{1495},\frac{2464}{5785}\right]$ & 2047,1495,5785 & 2464 & $\frac{3}{11}$ \Tstrut\Bstrut \\
$\left[\frac{1}{2},\frac{2640}{9401},\frac{2640}{7979},\frac{2640}{12019}\right]$ & 9401,7979,12019 & 2640 & $\frac{1}{10}$ \Tstrut\Bstrut \\
$\left[\frac{1}{2},\frac{3264}{20447},\frac{3264}{18415},\frac{3264}{23345}\right]$ & 20447,18415,23345 & 3264 & $\frac{1}{17}$ \Tstrut\Bstrut \\
$\left[\frac{1}{2},\frac{3520}{721},\frac{3520}{511},\frac{3520}{7519}\right]$ & 721,511,7519 & 3520 & $\frac{3}{8}$ \Tstrut\Bstrut \\
$\left[\frac{1}{2},\frac{3696}{4633},\frac{3696}{3485},\frac{3696}{9605}\right]$ & 4633,3485,9605 & 3696 & $\frac{2}{9}$ \Tstrut\Bstrut \\
$\left[\frac{1}{2},\frac{4160}{4879},\frac{4160}{3649},\frac{4160}{10591}\right]$ & 4879,3649,10591 & 4160 & $\frac{3}{13}$ \Tstrut\Bstrut \\
$\left[\frac{1}{2},\frac{6240}{959},\frac{6240}{679},\frac{6240}{13289}\right]$ & 959,679,13289 & 6240 & $\frac{5}{13}$ \Tstrut\Bstrut \\
$\left[\frac{1}{2},\frac{7280}{4681},\frac{7280}{3379},\frac{7280}{16459}\right]$ & 4681,3379,16459 & 7280 & $\frac{3}{10}$ \Tstrut\Bstrut \\
\bottomrule
\\[-0.5em]
\caption{Computer search results matching $s=\left[\frac{1}{2},\frac{d}{b},\frac{d}{a},\frac{d}{c}\right]$ }
\label{table:search}
\end{longtable}
\vspace{-2.5em}
In the table, $q$ is the variable in equation (\ref{eq:srat}). We also notice from the table that d is divisible by 16.

If we substitute $\left[\frac{1}{2},\frac{d}{b},\frac{d}{a},\frac{d}{c}\right]$ into equation (\ref{eq:govern}), we obtain
\begin{equation*}
\frac{d^4(2b^2c^2d^4-4a^2c^2d^4+2a^2b^2d^4-9a^2b^2c^2d^2+2a^2b^2c^4-4a^2b^4c^2+2a^4b^2c^2)}{8a^4b^4c^4} = 0
\end{equation*}
which is the relationship between $a,b,c,d$. We can solve for $d$ in terms of $a$,$b$,$c$, but an homogeneous sextic equation in $a$,$b$,$c$ results.

Without providing the derivation, the requirements for equation (\ref{eq:govern}) are met by the following two solution sets, except
that the reciprocals, $\frac{b}{d}$, $\frac{a}{d}$, $\frac{c}{d}$, are used instead, but the results are equivalent:

For $m,n \in \mathbb{Z}$:
\begin{equation}
s_{\mathbb{Z}} (m,n) = \left[\frac{1}{2}, \frac{-3m^4 + 18n^2m^2 - 3n^4}{8nm^3 - 8n^3m}, \frac{3m^4 + 6nm^3 + 6n^3m - 3n^4}{8nm^3 - 8n^3m}, \frac{3m^4 - 6nm^3 - 6n^3m - 3n^4}{8nm^3 - 8n^3m}\right]
\end{equation}

For $q \in \mathbb{Q}$:
\begin{equation}
s_{\mathbb{Q}} (q) = \left[\frac{1}{2},\frac{-3q^4+18q^2-3}{8q^3-8q},\frac{3q^4+6q^3+6q-3}{8q^3-8q},\frac{3q^4-6q^3-6q-3}{8q^3-8q}\right]
\end{equation}

Please note that from the first set, that $d = 8nm^3-8n^3m = 8mn(m-n)(m+n)$ is divisible by 16 as previously noted, because $mn(m-n)(m+n)$ is always divisible by 2,
for any $m,n \in \mathbb{Z}$.

\section*{{\bf A new parameterization for the $s$-parameter $\left[s_1,s_2,s_3,s_4\right]$ }}

After carefully examining over 2,000 $s$-parameters found by computer search, it was noticed that a particular $s$-parameter solution $\left[s_1,s_2,1,s_4\right]$
occurred from time to time. A parametric solution was found, which satisfied equation (\ref{eq:govern}), however none of these sets created a rational
Diophantine monoclinic piped, due to the fact that $s_3 = 1$ always leads to a degenerate solution.

After experimenting with these types of sets, a new solution for $s_1 = \frac{r}{s}$ was discovered.

A parametric solution in $\mathbb{Z}$:
\begin{align}
S_{\mathbb{Z}} (r,s,m,n) & = \left[s_1,s_2,s_3,s_4\right] \text{ for } r,s,m,n \in \mathbb{Z} \text{ where} \label{eq:sZ} \\
s_1 & = \frac{r}{s} \notag \\
s_2 & = \frac{(r^2-s^2)m^4+(-6r^2+6s^2)n^2m^2+(r^2-s^2)n^4}{4srnm^3-4srn^3m} \notag \\
s_3 & = \frac{(-r^2+s^2)m^4+(-2r^2+2s^2)nm^3+(-2r^2+2s^2)n^3m+(r^2-s^2)n^4}{4srnm^3-4srn^3m} \notag \\
s_4 & = \frac{(r^2-s^2)m^4+(-2r^2+2s^2)nm^3+(-2r^2+2s^2)n^3m+(-r^2+s^2)n^4}{4srn^3m-4srnm^3} \notag \\
T_{\mathbb{Z}} (r,s,m,n) & = \frac{(-r^2+s^2)m^2+(-r^2+s^2)n^2}{srm^2-srn^2} = s_3 - s_4
\end{align}
where $T_{\mathbb{Z}} (r,s,m,n) = s_3-s_4$ which is a measure of how close the monoclinic piped approaches a rectangular cuboid where $s_3 = s_4$.

A parametric solution in $\mathbb{Q}$:
\begin{align}
S_{\mathbb{Q}} (s,r) & = \left[s_1,s_2,s_3,s_4\right] \text{ for } s,r \in \mathbb{Q} \text{ where} \label{eq:sQ} \\
s_1 & = s \notag \\
s_2 & = \frac{(s^2-1)r^4+(-6s^2+6)r^2+(s^2-1)}{4sr^3-4sr} \notag \\
s_3 & = \frac{(-s^2+1)r^4+(-2s^2+2)r^3+(-2s^2+2)r+(s^2-1)}{4sr^3-4sr} \notag \\
s_4 & = \frac{(s^2-1)r^4+(-2s^2+2)r^3+(-2s^2+2)r+(-s^2+1)}{4sr-4sr^3} \notag \\
T_{\mathbb{Q}} (s,r) & = \frac{(-s^2+1)r^2+(-s^2+1)}{sr^2-s} = s_3 - s_4
\end{align}
and similarly, $T_{\mathbb{Q}} (s,r) = s_3 - s_4$, which is a measure of how close the monoclinic piped approaches
the rectangular cuboid.

\section*{{\bf Asymptotic Sequences for pipeds approaching the perfect cuboid }}

After discovering the parameterization, equation (\ref{eq:sQ}), which satisfies the second pattern type, the computer
was programmed to create 100's of millions of $s$-parameters sets for $S_{\mathbb{Q}}(s,r) = \left[s_1,s_2,s_3,s_4\right]$
for $s=\frac{m}{n}$ a proper fraction, $m<n<=500$, and 202,861 unique values of $r$.

A batch file was created, reduced, and sorted. 121,251,195 $s$-parameter sets were obtained.

It was decided to plot these $s$-parameter points, using the $x = s_3-s_4$ as the coordinate and $y = \text{numerator}(s_2)$
as the ordinate. The idea was to discover if perhaps the points had some type of pattern.

\begin{figure}[!ht]
\centering
\includegraphics[width=\textwidth]{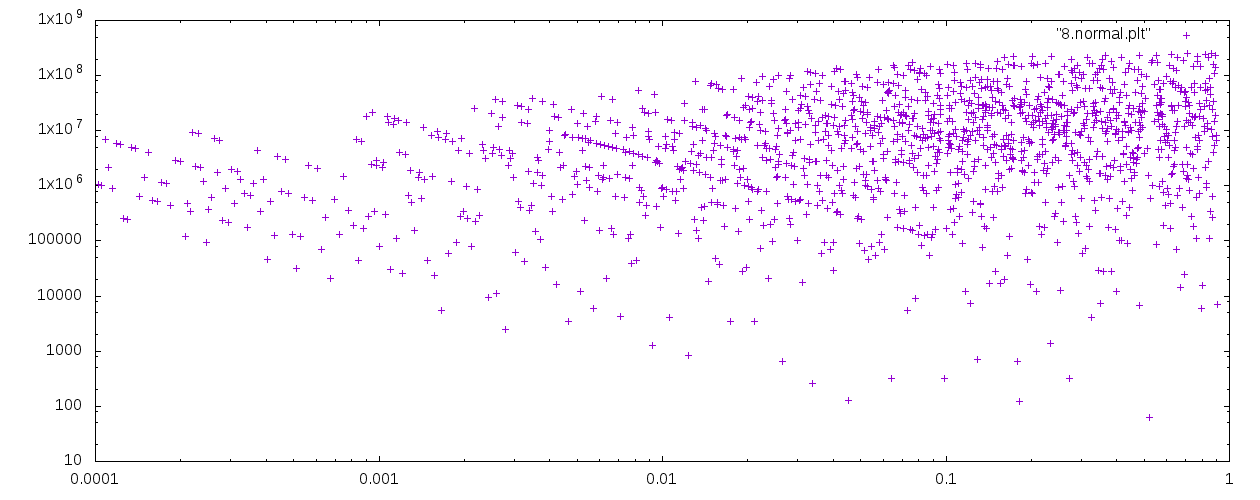}
\caption{Raw points from the parametric solution for $s_1 = \frac{1}{8}$}
\label{fig:oneeighthraw}
\end{figure}

Interestingly, after setting both the x-axis and y-axis to a log scale, the points do follow lines, strongly hinting that
they are found on rational polynomial fraction expressions. See Figure \ref{fig:oneeighthraw}.

\begin{figure}[!ht]
\centering
\includegraphics[width=\textwidth]{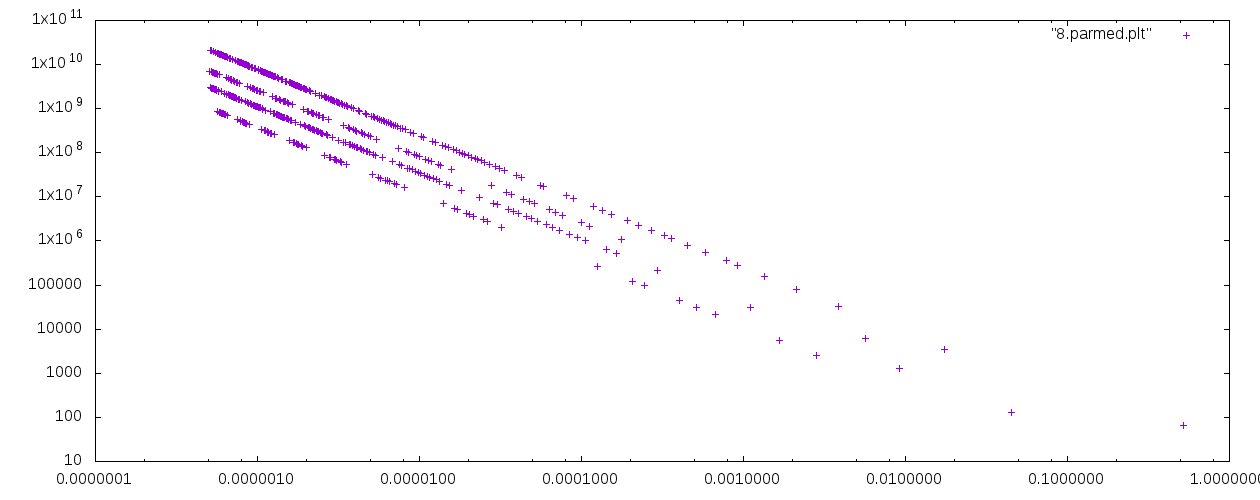}
\caption{Parameterized points from the parametric solution for $s_1 = \frac{1}{8}$}
\label{fig:oneeighthparm}
\end{figure}

This is indeed the case, as uncovering the solutions for $r=\frac{1}{8}$ shows in Figure \ref{fig:oneeighthparm}.

\begin{figure}[!ht]
\centering
\includegraphics[width=\textwidth]{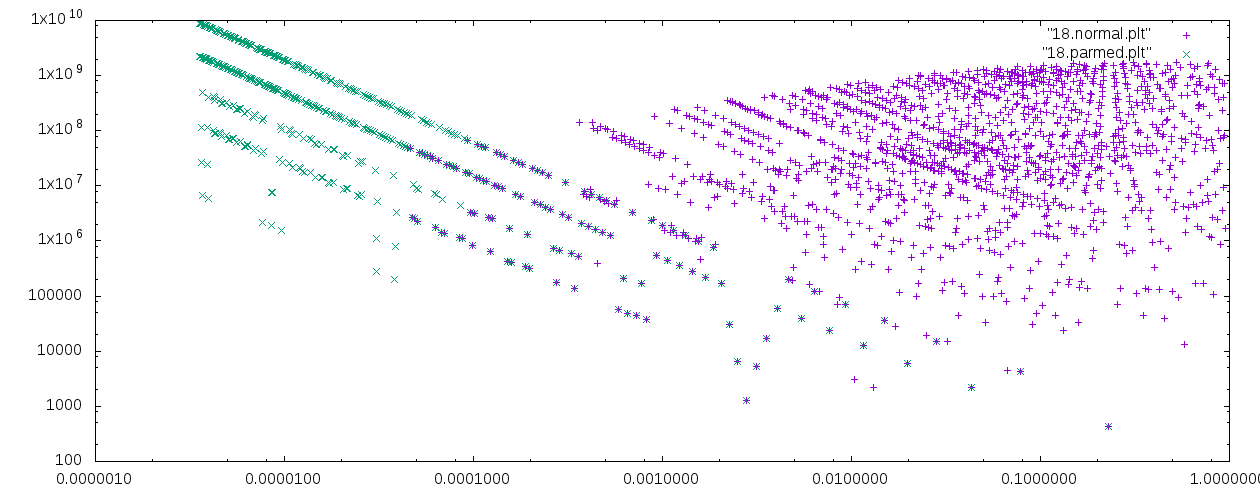}
\caption{Raw and parameterized points from the parametric solution for $s_1 = \frac{1}{18}$}
\label{fig:oneeighteenth}
\end{figure}

As another example, for $s_1 = \frac{1}{18}$, in Figure \ref{fig:oneeighteenth}, we can recognize what is happening.
These plots contain rational polynomial fraction expressions for the points on a line. Careful point fitting revealed
that the coordinates of the points were from rational polynomial fractions with a 3rd degree numerator and 4th degree denominator.

The plots enable the author to identify asymptotic sequences, for the monoclinic piped, where $s_3 \approxeq s_4$.
This means that the piped is approaching the {\it perfect cuboid}, either from an obtuse angle $> 90^\degree$ or
an acute angle $< 90^\degree$ as the angle $\to 90^\degree$.
\newpage
{\bf Obtuse angle asymptotic sequence}

The first discovered, because it was the easiest to recognize, was the obtuse angle asymptotic sequence where the
obtuse angle $\to 90^\degree$. In these equations, $n \in \mathbb{Z^+}$:
\begin{align}
s_{ob}(n) & = \left[\frac{n-1}{n},\frac{4n^4-8n^3+4n-1}{4n^4-8n^3+4n^2},\frac{4n^4-12n^3+12n^2-6n+1}{4n^4-8n^3+4n^2},\frac{4n^4-8n^3+4n^2}{4n^4-4n^3+2n-1}\right] \\
t_{ob}(n) & = \frac{-(2n-1)^4}{4(n-1)^2n^2(2n^2-1)(2n^2-2n+1)}
\end{align}

{\bf Acute angle asymptotic sequences}

The acute angle asymptotic sequence initially defied recognition, and it became necessary to gradually build up the bounding
sequences depending upon $s1 = \{ \frac{1}{2}, \frac{1}{3}, \frac{1}{4}, \frac{1}{5}, \dots \}$ in hopes of obtaining
the sequence.

These asymptotic sequences were eventually identified:

\begin{align}
s_{\frac{1}{2}}(n) & = \left[\frac{1}{2},\frac{8n(n+1)(n+2)}{3(n^2-2)(n^2+4n+2)},\frac{8n(n+1)(n+2)}{3(n^2-2)(n^2+2n+2)},\frac{8n(n+1)(n+2)}{3(n^2+2n+2)(n^2+4n+2)}\right] \\
t_{\frac{1}{2}}(n) & = \frac{32n(n+1)^2(n+2)}{3(n^2-2)(n^2+2n+2)(n^2+4n+2)}
\end{align}

\begin{align}
s_{\frac{1}{3}}(n) & = \left[\frac{1}{3},\frac{3n(n+1)(n+2)}{2(n^2-2)(n^2+4n+2)},\frac{3n(n+1)(n+2)}{2(n^2-2)(n^2+2n+2)},\frac{3n(n+1)(n+2)}{2(n^2+2n+2)(n^2+4n+2)}\right] \\
t_{\frac{1}{3}}(n) & = \frac{6n(n+1)^2(n+2)}{(n^2-2)(n^2+2n+2)(n^2+4n+2)}
\end{align}

\begin{align}
s_{\frac{1}{4}}(n) & = \left[\frac{1}{4},\frac{16n(n+1)(n+2)}{15(n^2-2)(n^2+4n+2)},\frac{16n(n+1)(n+2)}{15(n^2-2)(n^2+2n+2)},\frac{16n(n+1)(n+2)}{15(n^2+2n+2)(n^2+4n+2)}\right] \\
t_{\frac{1}{4}}(n) & = \frac{64n(n+1)^2(n+2)}{15(n^2-2)(n^2+2n+2)(n^2+4n+2)}
\end{align}

\begin{align}
s_{\frac{1}{5}}(n) & = \left[\frac{1}{5},\frac{5n(n+1)(n+2)}{6(n^2-2)(n^2+4n+2)},\frac{5n(n+1)(n+2)}{6(n^2-2)(n^2+2n+2)},\frac{5n(n+1)(n+2)}{6(n^2+2n+2)(n^2+4n+2)}\right] \\
t_{\frac{1}{5}}(n) & = \frac{10n(n+1)^2(n+2)}{3(n^2-2)(n^2+2n+2)(n^2+4n+2)}
\end{align}

\begin{align}
s_{\frac{1}{6}}(n) & = \left[\frac{1}{6},\frac{24n(n+1)(n+2)}{35(n^2-2)(n^2+4n+2)},\frac{24n(n+1)(n+2)}{35(n^2-2)(n^2+2n+2)},\frac{24n(n+1)(n+2)}{35(n^2+2n+2)(n^2+4n+2)}\right] \\
t_{\frac{1}{6}}(n) & = \frac{96n(n+1)^2(n+2)}{35(n^2-2)(n^2+2n+2)(n^2+4n+2)}
\end{align}

\begin{align}
s_{\frac{1}{7}}(n) & = \left[\frac{1}{7},\frac{7n(n+1)(n+2)}{12(n^2-2)(n^2+4n+2)},\frac{7n(n+1)(n+2)}{12(n^2-2)(n^2+2n+2)},\frac{7n(n+1)(n+2)}{12(n^2+2n+2)(n^2+4n+2)}\right] \\
t_{\frac{1}{7}}(n) & = \frac{7n(n+1)^2(n+2)}{3(n^2-2)(n^2+2n+2)(n^2+4n+2)}
\end{align}

\begin{align}
s_{\frac{1}{8}}(n) & = \left[\frac{1}{8},\frac{32n(n+1)(n+2)}{63(n^2-2)(n^2+4n+2)},\frac{32n(n+1)(n+2)}{63(n^2-2)(n^2+2n+2)},\frac{32n(n+1)(n+2)}{63(n^2+2n+2)(n^2+4n+2)}\right] \\
t_{\frac{1}{8}}(n) & = \frac{128n(n+1)^2(n+2)}{63(n^2-2)(n^2+2n+2)(n^2+4n+2)}
\end{align}

\begin{align}
s_{\frac{1}{9}}(n) & = \left[\frac{1}{9},\frac{9n(n+1)(n+2)}{20(n^2-2)(n^2+4n+2)},\frac{9n(n+1)(n+2)}{20(n^2-2)(n^2+2n+2)},\frac{9n(n+1)(n+2)}{20(n^2+2n+2)(n^2+4n+2)}\right] \\
t_{\frac{1}{9}}(n) & = \frac{9n(n+1)^2(n+2)}{5(n^2-2)(n^2+2n+2)(n^2+4n+2)}
\end{align}

\begin{align}
s_{\frac{1}{10}}(n) & = \left[\frac{1}{10},\frac{40n(n+1)(n+2)}{99(n^2-2)(n^2+4n+2)},\frac{40n(n+1)(n+2)}{99(n^2-2)(n^2+2n+2)},\frac{40n(n+1)(n+2)}{99(n^2+2n+2)(n^2+4n+2)}\right] \\
t_{\frac{1}{10}}(n) & = \frac{160n(n+1)^2(n+2)}{99(n^2-2)(n^2+2n+2)(n^2+4n+2)}
\end{align}

\begin{align}
s_{\frac{1}{11}}(n) & = \left[\frac{1}{11},\frac{11n(n+1)(n+2)}{30(n^2-2)(n^2+4n+2)},\frac{11n(n+1)(n+2)}{30(n^2-2)(n^2+2n+2)},\frac{11n(n+1)(n+2)}{30(n^2+2n+2)(n^2+4n+2)}\right] \\
t_{\frac{1}{11}}(n) & = \frac{22n(n+1)^2(n+2)}{15(n^2-2)(n^2+2n+2)(n^2+4n+2)}
\end{align}

\begin{align}
s_{\frac{1}{12}}(n) & = \left[\frac{1}{12},\frac{48n(n+1)(n+2)}{143(n^2-2)(n^2+4n+2)},\frac{48n(n+1)(n+2)}{143(n^2-2)(n^2+2n+2)},\frac{48n(n+1)(n+2)}{143(n^2+2n+2)(n^2+4n+2)}\right] \\
t_{\frac{1}{12}}(n) & = \frac{192n(n+1)^2(n+2)}{143(n^2-2)(n^2+2n+2)(n^2+4n+2)}
\end{align}

\begin{align}
s_{\frac{1}{13}}(n) & = \left[\frac{1}{13},\frac{13n(n+1)(n+2)}{42(n^2-2)(n^2+4n+2)},\frac{13n(n+1)(n+2)}{42(n^2-2)(n^2+2n+2)},\frac{13n(n+1)(n+2)}{42(n^2+2n+2)(n^2+4n+2)}\right] \\
t_{\frac{1}{13}}(n) & = \frac{26n(n+1)^2(n+2)}{21(n^2-2)(n^2+2n+2)(n^2+4n+2)}
\end{align}

\begin{align}
s_{\frac{1}{14}}(n) & = \left[\frac{1}{14},\frac{56n(n+1)(n+2)}{195(n^2-2)(n^2+4n+2)},\frac{56n(n+1)(n+2)}{195(n^2-2)(n^2+2n+2)},\frac{56n(n+1)(n+2)}{195(n^2+2n+2)(n^2+4n+2)}\right] \\
t_{\frac{1}{14}}(n) & = \frac{224n(n+1)^2(n+2)}{195(n^2-2)(n^2+2n+2)(n^2+4n+2)}
\end{align}

\begin{align}
s_{\frac{1}{15}}(n) & = \left[\frac{1}{15},\frac{15n(n+1)(n+2)}{56(n^2-2)(n^2+4n+2)},\frac{15n(n+1)(n+2)}{56(n^2-2)(n^2+2n+2)},\frac{15n(n+1)(n+2)}{56(n^2+2n+2)(n^2+4n+2)}\right] \\
t_{\frac{1}{15}}(n) & = \frac{15n(n+1)^2(n+2)}{14(n^2-2)(n^2+2n+2)(n^2+4n+2)}
\end{align}

\begin{align}
s_{\frac{1}{16}}(n) & = \left[\frac{1}{16},\frac{64n(n+1)(n+2)}{255(n^2-2)(n^2+4n+2)},\frac{64n(n+1)(n+2)}{255(n^2-2)(n^2+2n+2)},\frac{64n(n+1)(n+2)}{255(n^2+2n+2)(n^2+4n+2)}\right] \\
t_{\frac{1}{16}}(n) & = \frac{256n(n+1)^2(n+2)}{255(n^2-2)(n^2+2n+2)(n^2+4n+2)}
\end{align}

\begin{align}
s_{\frac{1}{17}}(n) & = \left[\frac{1}{17},\frac{17n(n+1)(n+2)}{72(n^2-2)(n^2+4n+2)},\frac{17n(n+1)(n+2)}{72(n^2-2)(n^2+2n+2)},\frac{17n(n+1)(n+2)}{72(n^2+2n+2)(n^2+4n+2)}\right] \\
t_{\frac{1}{17}}(n) & = \frac{17n(n+1)^2(n+2)}{18(n^2-2)(n^2+2n+2)(n^2+4n+2)}
\end{align}

\begin{align}
s_{\frac{1}{18}}(n) & = \left[\frac{1}{18},\frac{72n(n+1)(n+2)}{323(n^2-2)(n^2+4n+2)},\frac{72n(n+1)(n+2)}{323(n^2-2)(n^2+2n+2)},\frac{72n(n+1)(n+2)}{323(n^2+2n+2)(n^2+4n+2)}\right] \\
t_{\frac{1}{18}}(n) & = \frac{288n(n+1)^2(n+2)}{323(n^2-2)(n^2+2n+2)(n^2+4n+2)}
\end{align}

\begin{align}
s_{\frac{1}{19}}(n) & = \left[\frac{1}{19},\frac{19n(n+1)(n+2)}{90(n^2-2)(n^2+4n+2)},\frac{19n(n+1)(n+2)}{90(n^2-2)(n^2+2n+2)},\frac{19n(n+1)(n+2)}{90(n^2+2n+2)(n^2+4n+2)}\right] \\
t_{\frac{1}{19}}(n) & = \frac{38n(n+1)^2(n+2)}{45(n^2-2)(n^2+2n+2)(n^2+4n+2)}
\end{align}

After obtaining these 18 sequences, Maxima was used to derive the coefficients where $d\to\infty$.

{\bf A two parameter acute angle asymptotic sequence}

A solution for the asymptotic sequence where the acute angle $\to 90^\degree$ contains two
parameters for the input, $d,n \in \mathbb{Z^+}$, where $d$ is the value of the denominator in the $s_1 = \frac{1}{d}$ term,
and $n$ is the sequence number in that particular $s_1$ sequence.

This is not the same as the sequence for the obtuse angle $\to 90^\degree$ as it only takes 1 parameter for input.

The solution for the two parameter acute angle asymptotic sequence is:

\begin{align}
\begin{split}
s_{\frac{1}{d}}(d,n) & = \bigg [ \frac{1}{d}, \frac{4dn(n+1)(n+2)}{(d-1)(d+1)(n^2-2)(n^2+4n+2)}, \\
 & \qquad \quad \frac{4dn(n+1)(n+2)}{(d-1)(d+1)(n^2-2)(n^2+2n+2)}, \frac{4dn(n+1)(n+2)}{(d-1)(d+1)(n^2+2n+2)(n^2+4n+2)} \bigg ]
\end{split} \\
t_{\frac{1}{d}}(d,n) & = \frac{16dn(n+1)^2(n+2)}{(d-1)(d+1)(n^2-2)(n^2+2n+2)(n^2+4n+2)}
\end{align}
where $d,n \in \mathbb{Z^+}$.

\section*{{\bf Some concluding comments on the $s$-parameters and their pipeds }}

Initially the computer was programmed to do raw brute force searches for $s$-parameter sets for equation (\ref{eq:govern})
and found 3,280 solutions.

After discovering the parameter solution, equation (\ref{eq:sZ}) or equation(\ref{eq:sQ}), the computer created
a large batch of solutions to see how many of the raw solutions were included in this set.

It was discovered that only 16 solutions do not appear in the parameterization, thus the solution covered 99.5\% of
the raw solutions, which is quite remarkable.

These 16 $s$-parameter solutions are
\begin{longtable}{l}
$\left[ s_1, s_2, s_3, s_4 \right]$ \\
\midrule
$\left[\frac{3}{10},\frac{4}{15},\frac{4}{5},\frac{3}{20}\right]$ \Tstrut\Bstrut \\
$\left[\frac{42}{55},\frac{35}{132},\frac{20}{77},\frac{35}{132}\right]$ \Tstrut\Bstrut \\
$\left[\frac{41}{65},\frac{455}{943},\frac{13}{35},\frac{23}{41}\right]$ \Tstrut\Bstrut \\
$\left[\frac{143}{217},\frac{403}{616},\frac{8}{11},\frac{104}{217}\right]$ \Tstrut\Bstrut \\
$\left[\frac{68}{401},\frac{401}{721},\frac{2807}{4964},\frac{4964}{41303}\right]$ \Tstrut\Bstrut \\
$\left[\frac{348}{401},\frac{401}{527},\frac{6817}{8700},\frac{8700}{12431}\right]$ \Tstrut\Bstrut \\
$\left[\frac{314}{415},\frac{415}{527},\frac{1411}{1570},\frac{1570}{2573}\right]$ \Tstrut\Bstrut \\
$\left[\frac{188}{433},\frac{433}{623},\frac{964}{3031},\frac{38537}{45308}\right]$ \Tstrut\Bstrut \\
$\left[\frac{253}{439},\frac{439}{527},\frac{6325}{7463},\frac{6325}{13609}\right]$ \Tstrut\Bstrut \\
$\left[\frac{205}{457},\frac{457}{527},\frac{5125}{7769},\frac{5125}{14167}\right]$ \Tstrut\Bstrut \\
$\left[\frac{294}{473},\frac{473}{697},\frac{3738}{8041},\frac{19393}{26166}\right]$ \Tstrut\Bstrut \\
$\left[\frac{55}{479},\frac{479}{721},\frac{3353}{4015},\frac{4015}{49337}\right]$ \Tstrut\Bstrut \\
$\left[\frac{138}{481},\frac{481}{527},\frac{3450}{8177},\frac{3450}{14911}\right]$ \Tstrut\Bstrut \\
$\left[\frac{341}{661},\frac{661}{791},\frac{1903}{4627},\frac{58993}{74693}\right]$ \Tstrut\Bstrut \\
$\left[\frac{38}{751},\frac{751}{959},\frac{3686}{5257},\frac{3686}{102887}\right]$ \Tstrut\Bstrut \\
$\left[\frac{62}{769},\frac{769}{791},\frac{346}{5383},\frac{10726}{86897}\right]$ \Tstrut\Bstrut \\
\bottomrule
\\[-0.5em]
\caption{Raw search results not matching any parametrization.}
\label{table:anomalous}
\end{longtable}
\vspace{-1.0cm}

It would be interesting to come up with some new rational polynomial fraction expressions for these $s$-parameters,
as the 5th and following terms listed above seem to follow the same type of pattern, $s1,s2 = \frac{a}{b},\frac{b}{c}$.

As a general note, none of the rational polynomial fractions in this paper, for $s_3$, and $s_4$ can be equated in
an attempt to find the {\it perfect cuboid} solution, as doing so caused the parameter $\frac{r}{s}$ or $r$ to
go to 1, giving a degenerate solution. Thus no perfect cuboids exist as a solution of any of these expressions.

This author also checked all the monoclinic pipeds for rational area or volume, from the parameterizations, none had
any face parallelograms with rational area, nor did any have a rational volume. None of the Diophantine monoclinic
piped parametrizations have the vertices embedded in the rational lattice. This is in regard to two questions from
Sawyer and Reiter \cite{saw}.

\noindent \rule[0pt]{3.0in}{0.4pt}

\noindent \rule[0pt]{3.0in}{0.25pt}

\begin{thebibliography}{10}

	\bibitem{guy} Guy, Richard K., {\em Problem D18, Is there a perfect cuboid?},
	Unsolved Problems in Number Theory, 3rd ed. Springer-Verlag, NY.

	\bibitem{rath} Rathbun, Randall L., {\em Four integer parametrizations for the monoclinic Diophantine piped}, \\
	https://arXiv.org/abs/1705.07734 [math.NT] 2 June 2017

	\bibitem{reit1} Reiter, Clifford A., Tirrell, Jordan O., {\em Pursuing the Perfect Parallelepiped}, \\
	https://webbox.lafayette.edu/\~reiterc/nt/ppllpd/ppllpdpp.pdf

	\bibitem{reit2} Reiter, Clifford A., Tirrell, Jordan O. {\em Pursuing the Perfect Parallelepiped Auxiliary Materials}, \\
	http://webbox.lafayette.edu/\~reiterc/nt/ppllpd/index.html

	\bibitem{saw} Sawyer, Jorge F., Reiter, Clifford A., {\em Perfect parallelepipeds exist},
	Mathematics of Computation, vol 80, (2011), pp 1037-1040,
	also arXiv:0907.0220v2[math.NT]

	\bibitem{shar} Sharipov, Ruslan, {\em On Walter Wyss's no perfect cuboid paper}, \\
	https://arXiv.org/abs/1704.00165 [math.NT]

	\bibitem{soko} Sokolowshy, Benjamin D., VanHooft, Amy G., Volkert, Rachel M., Reiter, Clifford A.,
	{\em An Infinite Family of Perfect Parallelepipeds}, Mathematics of Computation, 83 (2014), pp 2441-2454.

	\bibitem{wyss1} Wyss, Walter, {\em Perfect Parallelograms},
	American Mathematical Monthly, 119 (6) (2012), pp 513-515

	\bibitem{wyss2} Wyss, Walter, {\em Sums of Squares, Bijective Parameter Representation}, \\
	https://arXiv.org/abs/1402.0102v1 [math.NT], 1 Feb 2014

	\bibitem{wyss3} Wyss, Walter, {\em No Perfect Cuboid},
	https://arXiv.org/abs/1506.02215v5 [math.NT], 17 May 2017

\end{thebibliography}
\end{document}